\documentclass[a4paper,12pt,reqno]{amsart}
\newtheorem{lem}{Lemma}
\newtheorem{theo}{Theorem}

\newtheorem{defi}{Definition}
\numberwithin{equation}{section}

\RequirePackage{graphics}

\def\t@xderror #1{%
   \GenericError{%
      \space\space\space\@spaces\@spaces\@spaces
   }{%
      TeXdraw Error: #1%
   }{%
      See the TeXdraw manual for an explanation.%
   }{\@ehc}%
}

\input texdraw
\newcommand{\RhombusVert}{\rlvec(0.8660254037844 0.5)
                          \rlvec(0.8660254037844 -0.5)
                          \rlvec(-0.8660254037844 -0.5)
                          \rlvec(-0.8660254037844  0.5)}
\newcommand{\RhombusPos}{\rlvec(0.8660254037844 0.5) 
                         \rlvec(0 -1)
                         \rlvec(-0.8660254037844 -0.5)
                         \rlvec(0 1)}   
\newcommand{\RhombusNeg}{\rlvec(0 1) 
                         \rlvec(0.8660254037844 -0.5)
                         \rlvec(0 -1)
                         \rlvec(-0.8660254037844 0.5)}  
\newcommand{\DreieckBreit}{\rlvec(0 1)
                           \rlvec(0.8660254037844 -0.5)
                           \rlvec(-0.8660254037844 -0.5)}
\newcommand{\DreieckSpitz}{\rlvec(0.8660254037844 0.5)
                           \rlvec(0 -1)
                           \rlvec(-0.8660254037844 0.5)}

\begin{document}

\title{Moments of inertia associated with the lozenge tilings of a hexagon}
\author{Ilse Fischer}

\address{Universit\"at Klagenfurt, Universit\"atsstrasse 65-67, A-9020 Klagenfurt, Austria.}

\email{ilse.fischer@uni-klu.ac.at}

\maketitle

\begin{abstract}
Consider the probability that an arbitrary chosen lozenge tiling of the 
hexagon with side lengths $a$, $b$, $c$, $a$, $b$, $c$ contains the 
horizontal lozenge with lowest vertex $(x,y)$ as if it described the 
distribution of mass in the plane. We compute the horizontal and 
the vertical moments of inertia with respect to this distribution. 
This solves a problem by Propp \cite[Problem 7]{propp}.
\end{abstract}

\section{Introduction}
\label{int}
Let $a$, $b$ and $c$ be positive integers and consider a hexagon with side lengths 
$a$,$b$,$c$,$a$,$b$,$c$ whose angles are $120^\circ$ (see Figure~\ref{fig:1}). 
The subject of our interest is lozenge tilings of such a
hexagon using lozenges with all sides of length $1$ and angles of $60^\circ$ and $120^\circ$. 
Figure~\ref{fig:2} shows an example of a lozenge tiling of a hexagon with $a=3$, $b=5$ 
and $c=4$. 

We introduce the following oblique angled coordinate system: Its origin is located 
in one of the two vertices, where sides of length $b$ and $c$ meet, and the axes
are induced by those two sides (see Figure~\ref{fig:3}). The units are chosen such that 
the side lengths of the considered hexagon are $\sqrt{2} a$, $b$, $c$, $\sqrt{2} a$, $b$, $c$ in this 
coordinate system. (That is to say, the two triangles in Figure~\ref{fig:3} with 
vertices in the origin form the unit `square'.)

\begin{figure}
\centertexdraw{
\drawdim cm
\DreieckBreit
\rmove(0 1)
\DreieckSpitz
\rmove(0.8660254037844 -0.5)
\DreieckBreit
\rmove(0 1)
\DreieckSpitz
\rmove(0.8660254037844 -0.5)
\DreieckBreit
\rmove(0 1)
\DreieckSpitz
\rmove(0.8660254037844 -0.5)
\DreieckBreit
\rmove(-0.8660254037844 -0.5)
\rmove(-0.8660254037844 -0.5)
\rmove(-0.8660254037844 -0.5)
\rmove(0 -1)
\DreieckBreit
\rmove(0 1)
\DreieckSpitz
\rmove(0.8660254037844 -0.5)
\DreieckBreit
\rmove(0 1)
\DreieckSpitz
\rmove(0.8660254037844 -0.5)
\DreieckBreit
\rmove(0 1)
\DreieckSpitz
\rmove(0.8660254037844 -0.5)
\DreieckBreit
\rmove(0 1)
\DreieckSpitz
\rmove(0.8660254037844 -0.5)
\DreieckBreit
\rmove(-0.8660254037844 -0.5)
\rmove(-0.8660254037844 -0.5)
\rmove(-0.8660254037844 -0.5)
\rmove(-0.8660254037844 -0.5)
\rmove(0 -1)
\DreieckBreit
\rmove(0 1)
\DreieckSpitz
\rmove(0.8660254037844 -0.5)
\DreieckBreit
\rmove(0 1)
\DreieckSpitz
\rmove(0.8660254037844 -0.5)
\DreieckBreit
\rmove(0 1)
\DreieckSpitz
\rmove(0.8660254037844 -0.5)
\DreieckBreit
\rmove(0 1)
\DreieckSpitz
\rmove(0.8660254037844 -0.5)
\DreieckBreit
\rmove(0 1)
\DreieckSpitz
\rmove(0.8660254037844 -0.5)
\DreieckBreit
\rmove(-0.8660254037844 -0.5)
\rmove(-0.8660254037844 -0.5)
\rmove(-0.8660254037844 -0.5)
\rmove(-0.8660254037844 -0.5)
\rmove(-0.8660254037844 -0.5)
\rmove(0 -1)
\DreieckBreit
\rmove(0 1)
\DreieckSpitz
\rmove(0.8660254037844 -0.5)
\DreieckBreit
\rmove(0 1)
\DreieckSpitz
\rmove(0.8660254037844 -0.5)
\DreieckBreit
\rmove(0 1)
\DreieckSpitz
\rmove(0.8660254037844 -0.5)
\DreieckBreit
\rmove(0 1)
\DreieckSpitz
\rmove(0.8660254037844 -0.5)
\DreieckBreit
\rmove(0 1)
\DreieckSpitz
\rmove(0.8660254037844 -0.5)
\DreieckBreit
\rmove(0 1)
\DreieckSpitz
\rmove(0.8660254037844 -0.5)
\DreieckBreit
\rmove(-0.8660254037844 -0.5)
\rmove(-0.8660254037844 -0.5)
\rmove(-0.8660254037844 -0.5)
\rmove(-0.8660254037844 -0.5)
\rmove(-0.8660254037844 -0.5)
\rmove(-0.8660254037844 -0.5)
\DreieckSpitz
\rmove(0.8660254037844 -0.5)
\DreieckBreit
\rmove(0 1)
\DreieckSpitz
\rmove(0.8660254037844 -0.5)
\DreieckBreit
\rmove(0 1)
\DreieckSpitz
\rmove(0.8660254037844 -0.5)
\DreieckBreit
\rmove(0 1)
\DreieckSpitz
\rmove(0.8660254037844 -0.5)
\DreieckBreit
\rmove(0 1)
\DreieckSpitz
\rmove(0.8660254037844 -0.5)
\DreieckBreit
\rmove(0 1)
\DreieckSpitz
\rmove(0.8660254037844 -0.5)
\DreieckBreit
\rmove(0 1)
\DreieckSpitz
\rmove(0.8660254037844 -0.5)
\DreieckBreit
\rmove(-0.8660254037844 -0.5)
\rmove(-0.8660254037844 -0.5)
\rmove(-0.8660254037844 -0.5)
\rmove(-0.8660254037844 -0.5)
\rmove(-0.8660254037844 -0.5)
\rmove(-0.8660254037844 -0.5)
\DreieckSpitz
\rmove(0.8660254037844 -0.5)
\DreieckBreit
\rmove(0 1)
\DreieckSpitz
\rmove(0.8660254037844 -0.5)
\DreieckBreit
\rmove(0 1)
\DreieckSpitz
\rmove(0.8660254037844 -0.5)
\DreieckBreit
\rmove(0 1)
\DreieckSpitz
\rmove(0.8660254037844 -0.5)
\DreieckBreit
\rmove(0 1)
\DreieckSpitz
\rmove(0.8660254037844 -0.5)
\DreieckBreit
\rmove(0 1)
\DreieckSpitz
\rmove(0.8660254037844 -0.5)
\DreieckBreit
\rmove(0 1)
\DreieckSpitz
\rmove(0.8660254037844 -0.5)
\rmove(-0.8660254037844 -0.5)
\rmove(-0.8660254037844 -0.5)
\rmove(-0.8660254037844 -0.5)
\rmove(-0.8660254037844 -0.5)
\rmove(-0.8660254037844 -0.5)
\rmove(-0.8660254037844 -0.5)
\DreieckSpitz
\rmove(0.8660254037844 -0.5)
\DreieckBreit
\rmove(0 1)
\DreieckSpitz
\rmove(0.8660254037844 -0.5)
\DreieckBreit
\rmove(0 1)
\DreieckSpitz
\rmove(0.8660254037844 -0.5)
\DreieckBreit
\rmove(0 1)
\DreieckSpitz
\rmove(0.8660254037844 -0.5)
\DreieckBreit
\rmove(0 1)
\DreieckSpitz
\rmove(0.8660254037844 -0.5)
\DreieckBreit
\rmove(0 1)
\DreieckSpitz
\rmove(0.8660254037844 -0.5)
\rmove(-0.8660254037844 -0.5)
\rmove(-0.8660254037844 -0.5)
\rmove(-0.8660254037844 -0.5)
\rmove(-0.8660254037844 -0.5)
\rmove(-0.8660254037844 -0.5)
\DreieckSpitz
\rmove(0.8660254037844 -0.5)
\DreieckBreit
\rmove(0 1)
\DreieckSpitz
\rmove(0.8660254037844 -0.5)
\DreieckBreit
\rmove(0 1)
\DreieckSpitz
\rmove(0.8660254037844 -0.5)
\DreieckBreit
\rmove(0 1)
\DreieckSpitz
\rmove(0.8660254037844 -0.5)
\DreieckBreit
\rmove(0 1)
\DreieckSpitz
\rmove(0.8660254037844 -0.5)
\rmove(-0.8660254037844 -0.5)
\rmove(-0.8660254037844 -0.5)
\rmove(-0.8660254037844 -0.5)
\rmove(-0.8660254037844 -0.5)
\DreieckSpitz
\rmove(0.8660254037844 -0.5)
\DreieckBreit
\rmove(0 1)
\DreieckSpitz
\rmove(0.8660254037844 -0.5)
\DreieckBreit
\rmove(0 1)
\DreieckSpitz
\rmove(0.8660254037844 -0.5)
\DreieckBreit
\rmove(0 1)
\DreieckSpitz
\rmove(0.8660254037844 -0.5)
\move(-0.5 -1)
\htext{$c$}
\move(1.299038105677 1.75)
\rmove(0  0.5)
\htext{$a$}
\rmove(0 -0.5)
\rmove(1.299038105677 0.75)
\rmove(1.299038105677 -0.75)
\rmove(0.8660254037844 -0.5)
\rmove(0 0.5)
\htext{$b$}
\rmove(0 -0.5)
\rmove(1.299038105677 -0.75)
\rmove(0.8660254037844 -0.5)
\rmove(0 -2)
\rmove(0.5 0)
\htext{$c$}
\rmove(-0.5 0)
\rmove(0 -2)
\rmove(-1.299038105677 -0.75)
\rmove(0 -0.5)
\htext{$a$}
\rmove(0.5 0)
\rmove(-1.299038105677 -0.75)
\rmove(-1.299038105677 0.75)
\rmove(-0.8660254037844 0.5)
\rmove(-0.8660254037844 0.5)
\rmove(0 -0.5)
\htext{$b$}
}
\caption{}
\label{fig:1}
\end{figure}

\begin{figure}
\centertexdraw{
\drawdim cm
\linewd 0.08
\rmove(-5 0)
\RhombusPos
\rmove(0.8660254037844 -0.5)
\RhombusNeg
\rmove(0 1)
\RhombusVert
\rmove(0.8660254037844 0.5)
\RhombusVert

\rmove(-0.8660254037844 -0.5)
\rmove(-0.8660254037844 -0.5)
\rmove(0 -1)

\RhombusPos
\rmove(0.8660254037844 -0.5)
\RhombusNeg
\rmove(0.8660254037844 0.5)
\rmove(0 1)
\RhombusPos
\rmove(0.8660254037844 0.5)
\RhombusPos

\rmove(-0.8660254037844 -0.5)
\rmove(-0.8660254037844 -0.5)
\rmove(-0.8660254037844 -0.5)
\rmove(0 -1)

\RhombusPos
\rmove(0.8660254037844 -0.5)
\RhombusNeg
\rmove(0.8660254037844 0.5)
\RhombusNeg
\rmove(0 1)
\RhombusVert
\rmove(0.8660254037844 0.5)
\RhombusVert
\rmove(0.8660254037844 0.5)
\RhombusNeg

\rmove(-0.8660254037844 -0.5)
\rmove(-0.8660254037844 -0.5)
\rmove(-0.8660254037844 -0.5)
\rmove(-0.8660254037844 -0.5)
\rmove(0 -1)

\RhombusVert
\rmove(0.8660254037844 0.5)
\rmove(0.8660254037844 -0.5)
\RhombusNeg
\rmove(0.8660254037844 0.5)
\rmove(0 1)
\RhombusPos
\rmove(0.8660254037844 -0.5)
\RhombusNeg
\rmove(0 1)
\RhombusVert
\rmove(0.8660254037844 0.5)
\RhombusNeg

\rmove(-0.8660254037844 -0.5)
\rmove(-0.8660254037844 -0.5)
\rmove(-0.8660254037844 -0.5)
\rmove(-0.8660254037844 -0.5)
\rmove(-0.8660254037844 -0.5)
\rmove(0 -1)

\RhombusNeg
\rmove(0.8660254037844 0.5)
\RhombusPos
\rmove(0.8660254037844 -0.5)
\RhombusNeg
\rmove(0.8660254037844 0.5)
\RhombusNeg
\rmove(0 1)
\RhombusVert
\rmove(0.8660254037844 0.5)
\rmove(0.8660254037844 0.5)
\RhombusPos
\rmove(0.8660254037844 0.5)
\RhombusNeg
\rmove(0.8660254037844 -0.5)
\RhombusNeg

\rmove(-0.8660254037844 -0.5)
\rmove(-0.8660254037844 -0.5)
\rmove(-0.8660254037844 -0.5)
\rmove(-0.8660254037844 -0.5)
\rmove(-0.8660254037844 -0.5)
\rmove(-0.8660254037844 -0.5)

\RhombusVert
\rmove(0.8660254037844 0.5)
\rmove(0.8660254037844 -0.5)
\RhombusNeg
\rmove(0.8660254037844 0.5)
\RhombusNeg
\rmove(0 1)
\RhombusVert
\rmove(0.8660254037844 0.5)
\RhombusVert
\rmove(0.8660254037844 0.5)
\RhombusNeg
\rmove(0.8660254037844 -0.5)
\RhombusNeg

\rmove(-0.8660254037844 -0.5)
\rmove(-0.8660254037844 -0.5)
\rmove(-0.8660254037844 -0.5)
\rmove(-0.8660254037844 -0.5)
\rmove(-0.8660254037844 -0.5)

\RhombusVert
\rmove(0.8660254037844 0.5)
\rmove(0.8660254037844 -0.5)
\RhombusNeg
\rmove(0.8660254037844 0.5)
\rmove(0 1)
\RhombusPos
\rmove(0.8660254037844 0.5)
\RhombusPos
\rmove(0.8660254037844 -0.5)
\RhombusNeg

\rmove(-0.8660254037844 -0.5)
\rmove(-0.8660254037844 -0.5)
\rmove(-0.8660254037844 -0.5)
\rmove(-0.8660254037844 -0.5)

\RhombusVert 
\rmove(0.8660254037844 0.5)
\rmove(0.8660254037844 0.5)
\RhombusPos
\rmove(0.8660254037844 0.5)
\RhombusVert

\rmove(-0.8660254037844 -0.5)
\rmove(-0.8660254037844 -0.5)
\rmove(0 -1)

\RhombusVert
\rmove(0.8660254037844 0.5)
\RhombusVert
\rmove(0.8660254037844 0.5)
\RhombusNeg
\rmove(0.8660254037844 0.5)
\RhombusPos

\move(0 -1)
\rmove(-5 0)
\linewd 0.01
\DreieckBreit
\rmove(0 1)
\DreieckSpitz
\rmove(0.8660254037844 -0.5)
\DreieckBreit
\rmove(0 1)
\DreieckSpitz
\rmove(0.8660254037844 -0.5)
\DreieckBreit
\rmove(0 1)
\DreieckSpitz
\rmove(0.8660254037844 -0.5)
\DreieckBreit

\rmove(-0.8660254037844 -0.5)
\rmove(-0.8660254037844 -0.5)
\rmove(-0.8660254037844 -0.5)
\rmove(0 -1)

\DreieckBreit
\rmove(0 1)
\DreieckSpitz
\rmove(0.8660254037844 -0.5)
\DreieckBreit
\rmove(0 1)
\DreieckSpitz
\rmove(0.8660254037844 -0.5)
\DreieckBreit
\rmove(0 1)
\DreieckSpitz
\rmove(0.8660254037844 -0.5)
\DreieckBreit
\rmove(0 1)
\DreieckSpitz
\rmove(0.8660254037844 -0.5)
\DreieckBreit

\rmove(-0.8660254037844 -0.5)
\rmove(-0.8660254037844 -0.5)
\rmove(-0.8660254037844 -0.5)
\rmove(-0.8660254037844 -0.5)
\rmove(0 -1)

\DreieckBreit
\rmove(0 1)
\DreieckSpitz
\rmove(0.8660254037844 -0.5)
\DreieckBreit
\rmove(0 1)
\DreieckSpitz
\rmove(0.8660254037844 -0.5)
\DreieckBreit
\rmove(0 1)
\DreieckSpitz
\rmove(0.8660254037844 -0.5)
\DreieckBreit
\rmove(0 1)
\DreieckSpitz
\rmove(0.8660254037844 -0.5)
\DreieckBreit
\rmove(0 1)
\DreieckSpitz
\rmove(0.8660254037844 -0.5)
\DreieckBreit

\rmove(-0.8660254037844 -0.5)
\rmove(-0.8660254037844 -0.5)
\rmove(-0.8660254037844 -0.5)
\rmove(-0.8660254037844 -0.5)
\rmove(-0.8660254037844 -0.5)
\rmove(0 -1)

\DreieckBreit
\rmove(0 1)
\DreieckSpitz
\rmove(0.8660254037844 -0.5)
\DreieckBreit
\rmove(0 1)
\DreieckSpitz
\rmove(0.8660254037844 -0.5)
\DreieckBreit
\rmove(0 1)
\DreieckSpitz
\rmove(0.8660254037844 -0.5)
\DreieckBreit
\rmove(0 1)
\DreieckSpitz
\rmove(0.8660254037844 -0.5)
\DreieckBreit
\rmove(0 1)
\DreieckSpitz
\rmove(0.8660254037844 -0.5)
\DreieckBreit
\rmove(0 1)
\DreieckSpitz
\rmove(0.8660254037844 -0.5)
\DreieckBreit

\rmove(-0.8660254037844 -0.5)
\rmove(-0.8660254037844 -0.5)
\rmove(-0.8660254037844 -0.5)
\rmove(-0.8660254037844 -0.5)
\rmove(-0.8660254037844 -0.5)
\rmove(-0.8660254037844 -0.5)

\DreieckSpitz
\rmove(0.8660254037844 -0.5)
\DreieckBreit
\rmove(0 1)
\DreieckSpitz
\rmove(0.8660254037844 -0.5)
\DreieckBreit
\rmove(0 1)
\DreieckSpitz
\rmove(0.8660254037844 -0.5)
\DreieckBreit
\rmove(0 1)
\DreieckSpitz
\rmove(0.8660254037844 -0.5)
\DreieckBreit
\rmove(0 1)
\DreieckSpitz
\rmove(0.8660254037844 -0.5)
\DreieckBreit
\rmove(0 1)
\DreieckSpitz
\rmove(0.8660254037844 -0.5)
\DreieckBreit
\rmove(0 1)
\DreieckSpitz
\rmove(0.8660254037844 -0.5)
\DreieckBreit

\rmove(-0.8660254037844 -0.5)
\rmove(-0.8660254037844 -0.5)
\rmove(-0.8660254037844 -0.5)
\rmove(-0.8660254037844 -0.5)
\rmove(-0.8660254037844 -0.5)
\rmove(-0.8660254037844 -0.5)

\DreieckSpitz
\rmove(0.8660254037844 -0.5)
\DreieckBreit
\rmove(0 1)
\DreieckSpitz
\rmove(0.8660254037844 -0.5)
\DreieckBreit
\rmove(0 1)
\DreieckSpitz
\rmove(0.8660254037844 -0.5)
\DreieckBreit
\rmove(0 1)
\DreieckSpitz
\rmove(0.8660254037844 -0.5)
\DreieckBreit
\rmove(0 1)
\DreieckSpitz
\rmove(0.8660254037844 -0.5)
\DreieckBreit
\rmove(0 1)
\DreieckSpitz
\rmove(0.8660254037844 -0.5)
\DreieckBreit
\rmove(0 1)
\DreieckSpitz
\rmove(0.8660254037844 -0.5)

\rmove(-0.8660254037844 -0.5)
\rmove(-0.8660254037844 -0.5)
\rmove(-0.8660254037844 -0.5)
\rmove(-0.8660254037844 -0.5)
\rmove(-0.8660254037844 -0.5)
\rmove(-0.8660254037844 -0.5)

\DreieckSpitz
\rmove(0.8660254037844 -0.5)
\DreieckBreit
\rmove(0 1)
\DreieckSpitz
\rmove(0.8660254037844 -0.5)
\DreieckBreit
\rmove(0 1)
\DreieckSpitz
\rmove(0.8660254037844 -0.5)
\DreieckBreit
\rmove(0 1)
\DreieckSpitz
\rmove(0.8660254037844 -0.5)
\DreieckBreit
\rmove(0 1)
\DreieckSpitz
\rmove(0.8660254037844 -0.5)
\DreieckBreit
\rmove(0 1)
\DreieckSpitz
\rmove(0.8660254037844 -0.5)

\rmove(-0.8660254037844 -0.5)
\rmove(-0.8660254037844 -0.5)
\rmove(-0.8660254037844 -0.5)
\rmove(-0.8660254037844 -0.5)
\rmove(-0.8660254037844 -0.5)

\DreieckSpitz
\rmove(0.8660254037844 -0.5)
\DreieckBreit
\rmove(0 1)
\DreieckSpitz
\rmove(0.8660254037844 -0.5)
\DreieckBreit
\rmove(0 1)
\DreieckSpitz
\rmove(0.8660254037844 -0.5)
\DreieckBreit
\rmove(0 1)
\DreieckSpitz
\rmove(0.8660254037844 -0.5)
\DreieckBreit
\rmove(0 1)
\DreieckSpitz
\rmove(0.8660254037844 -0.5)

\rmove(-0.8660254037844 -0.5)
\rmove(-0.8660254037844 -0.5)
\rmove(-0.8660254037844 -0.5)
\rmove(-0.8660254037844 -0.5)

\DreieckSpitz
\rmove(0.8660254037844 -0.5)
\DreieckBreit
\rmove(0 1)
\DreieckSpitz
\rmove(0.8660254037844 -0.5)
\DreieckBreit
\rmove(0 1)
\DreieckSpitz
\rmove(0.8660254037844 -0.5)
\DreieckBreit
\rmove(0 1)
\DreieckSpitz
\rmove(0.8660254037844 -0.5)
}
\caption{}
\label{fig:2}
\end{figure}

\begin{figure}
\centertexdraw{
\drawdim cm
\arrowheadtype t:F
\linewd 0.05
\move(0 -3)
\ravec(0 5)
\move(0 -3)
\rlvec(0.8660254037844 -0.5)
\rlvec(0 1)
\rlvec(-0.866025403784 0.5)
\rlvec(0 -1)
\rlvec(0.8660254037844 -0.5)
\rlvec(0.8660254037844 -0.5)
\rlvec(0.8660254037844 -0.5)
\rlvec(0.8660254037844 -0.5)
\rlvec(0.8660254037844 -0.5)
\ravec(0.8660254037844 -0.5)
\linewd 0.01
\move(0 0)
\DreieckBreit
\rmove(0 1)
\DreieckSpitz
\rmove(0.8660254037844 -0.5)
\DreieckBreit
\rmove(0 1)
\DreieckSpitz
\rmove(0.8660254037844 -0.5)
\DreieckBreit
\rmove(0 1)
\DreieckSpitz
\rmove(0.8660254037844 -0.5)
\DreieckBreit
\rmove(-0.8660254037844 -0.5)
\rmove(-0.8660254037844 -0.5)
\rmove(-0.8660254037844 -0.5)
\rmove(0 -1)
\DreieckBreit
\rmove(0 1)
\DreieckSpitz
\rmove(0.8660254037844 -0.5)
\DreieckBreit
\rmove(0 1)
\DreieckSpitz
\rmove(0.8660254037844 -0.5)
\DreieckBreit
\rmove(0 1)
\DreieckSpitz
\rmove(0.8660254037844 -0.5)
\DreieckBreit
\rmove(0 1)
\DreieckSpitz
\rmove(0.8660254037844 -0.5)
\DreieckBreit
\rmove(-0.8660254037844 -0.5)
\rmove(-0.8660254037844 -0.5)
\rmove(-0.8660254037844 -0.5)
\rmove(-0.8660254037844 -0.5)
\rmove(0 -1)
\DreieckBreit
\rmove(0 1)
\DreieckSpitz
\rmove(0.8660254037844 -0.5)
\DreieckBreit
\rmove(0 1)
\DreieckSpitz
\rmove(0.8660254037844 -0.5)
\DreieckBreit
\rmove(0 1)
\DreieckSpitz
\rmove(0.8660254037844 -0.5)
\DreieckBreit
\rmove(0 1)
\DreieckSpitz
\rmove(0.8660254037844 -0.5)
\DreieckBreit
\rmove(0 1)
\DreieckSpitz
\rmove(0.8660254037844 -0.5)
\DreieckBreit
\rmove(-0.8660254037844 -0.5)
\rmove(-0.8660254037844 -0.5)
\rmove(-0.8660254037844 -0.5)
\rmove(-0.8660254037844 -0.5)
\rmove(-0.8660254037844 -0.5)
\rmove(0 -1)
\DreieckBreit
\rmove(0 1)
\DreieckSpitz
\rmove(0.8660254037844 -0.5)
\DreieckBreit
\rmove(0 1)
\DreieckSpitz
\rmove(0.8660254037844 -0.5)
\DreieckBreit
\rmove(0 1)
\DreieckSpitz
\rmove(0.8660254037844 -0.5)
\DreieckBreit
\rmove(0 1)
\DreieckSpitz
\rmove(0.8660254037844 -0.5)
\DreieckBreit
\rmove(0 1)
\DreieckSpitz
\rmove(0.8660254037844 -0.5)
\DreieckBreit
\rmove(0 1)
\DreieckSpitz
\rmove(0.8660254037844 -0.5)
\DreieckBreit
\rmove(-0.8660254037844 -0.5)
\rmove(-0.8660254037844 -0.5)
\rmove(-0.8660254037844 -0.5)
\rmove(-0.8660254037844 -0.5)
\rmove(-0.8660254037844 -0.5)
\rmove(-0.8660254037844 -0.5)
\DreieckSpitz
\rmove(0.8660254037844 -0.5)
\DreieckBreit
\rmove(0 1)
\DreieckSpitz
\rmove(0.8660254037844 -0.5)
\DreieckBreit
\rmove(0 1)
\DreieckSpitz
\rmove(0.8660254037844 -0.5)
\DreieckBreit
\rmove(0 1)
\DreieckSpitz
\rmove(0.8660254037844 -0.5)
\DreieckBreit
\rmove(0 1)
\DreieckSpitz
\rmove(0.8660254037844 -0.5)
\DreieckBreit
\rmove(0 1)
\DreieckSpitz
\rmove(0.8660254037844 -0.5)
\DreieckBreit
\rmove(0 1)
\DreieckSpitz
\rmove(0.8660254037844 -0.5)
\DreieckBreit
\rmove(-0.8660254037844 -0.5)
\rmove(-0.8660254037844 -0.5)
\rmove(-0.8660254037844 -0.5)
\rmove(-0.8660254037844 -0.5)
\rmove(-0.8660254037844 -0.5)
\rmove(-0.8660254037844 -0.5)
\DreieckSpitz
\rmove(0.8660254037844 -0.5)
\DreieckBreit
\rmove(0 1)
\DreieckSpitz
\rmove(0.8660254037844 -0.5)
\DreieckBreit
\rmove(0 1)
\DreieckSpitz
\rmove(0.8660254037844 -0.5)
\DreieckBreit
\rmove(0 1)
\DreieckSpitz
\rmove(0.8660254037844 -0.5)
\DreieckBreit
\rmove(0 1)
\DreieckSpitz
\rmove(0.8660254037844 -0.5)
\DreieckBreit
\rmove(0 1)
\DreieckSpitz
\rmove(0.8660254037844 -0.5)
\DreieckBreit
\rmove(0 1)
\DreieckSpitz
\rmove(0.8660254037844 -0.5)
\rmove(-0.8660254037844 -0.5)
\rmove(-0.8660254037844 -0.5)
\rmove(-0.8660254037844 -0.5)
\rmove(-0.8660254037844 -0.5)
\rmove(-0.8660254037844 -0.5)
\rmove(-0.8660254037844 -0.5)
\DreieckSpitz
\rmove(0.8660254037844 -0.5)
\DreieckBreit
\rmove(0 1)
\DreieckSpitz
\rmove(0.8660254037844 -0.5)
\DreieckBreit
\rmove(0 1)
\DreieckSpitz
\rmove(0.8660254037844 -0.5)
\DreieckBreit
\rmove(0 1)
\DreieckSpitz
\rmove(0.8660254037844 -0.5)
\DreieckBreit
\rmove(0 1)
\DreieckSpitz
\rmove(0.8660254037844 -0.5)
\DreieckBreit
\rmove(0 1)
\DreieckSpitz
\rmove(0.8660254037844 -0.5)
\rmove(-0.8660254037844 -0.5)
\rmove(-0.8660254037844 -0.5)
\rmove(-0.8660254037844 -0.5)
\rmove(-0.8660254037844 -0.5)
\rmove(-0.8660254037844 -0.5)
\DreieckSpitz
\rmove(0.8660254037844 -0.5)
\DreieckBreit
\rmove(0 1)
\DreieckSpitz
\rmove(0.8660254037844 -0.5)
\DreieckBreit
\rmove(0 1)
\DreieckSpitz
\rmove(0.8660254037844 -0.5)
\DreieckBreit
\rmove(0 1)
\DreieckSpitz
\rmove(0.8660254037844 -0.5)
\DreieckBreit
\rmove(0 1)
\DreieckSpitz
\rmove(0.8660254037844 -0.5)
\rmove(-0.8660254037844 -0.5)
\rmove(-0.8660254037844 -0.5)
\rmove(-0.8660254037844 -0.5)
\rmove(-0.8660254037844 -0.5)
\DreieckSpitz
\rmove(0.8660254037844 -0.5)
\DreieckBreit
\rmove(0 1)
\DreieckSpitz
\rmove(0.8660254037844 -0.5)
\DreieckBreit
\rmove(0 1)
\DreieckSpitz
\rmove(0.8660254037844 -0.5)
\DreieckBreit
\rmove(0 1)
\DreieckSpitz
\rmove(0.8660254037844 -0.5)
\move(-0.5 -1)
\htext{$c$}
\move(1.299038105677 1.75)
\rmove(0  0.5)
\htext{$a$}
\rmove(0 -0.5)
\rmove(1.299038105677 0.75)
\rmove(1.299038105677 -0.75)
\rmove(0.8660254037844 -0.5)
\rmove(0 0.5)
\htext{$b$}
\rmove(0 -0.5)
\rmove(1.299038105677 -0.75)
\rmove(0.8660254037844 -0.5)
\rmove(0 -2)
\rmove(0.5 0)
\htext{$c$}
\rmove(-0.5 0)
\rmove(0 -2)
\rmove(-1.299038105677 -0.75)
\rmove(0 -0.5)
\htext{$a$}
\rmove(0.5 0)
\rmove(-1.299038105677 -0.75)
\rmove(-1.299038105677 0.75)
\rmove(-0.8660254037844 0.5)
\rmove(-0.8660254037844 0.5)
\rmove(0 -0.5)
\htext{$b$}
}
\caption{}
\label{fig:3}
\end{figure} 

Let $P_{a,b,c}(x,y)$ denote the probability that 
an arbitrary chosen lozenge tiling of the hexagon with side lengths $a$,$b$,$c$,$a$,$b$,$c$ contains
the horizontal lozenge with lowest vertex $(x,y)$ in the oblique angled coordinate
system. Note that $S=((a+b)/2, (a+c)/2)$ is the centre of the 
hexagon in question. Consider the probability $P_{a,b,c}(x,y)$ as if it 
described the distribution of mass. In \cite[Problem 7]{propp} Propp suggests to compute 
the horizontal moment of inertia with respect to $S$
$$
\sum_{y=0}^{a+c-1} \sum_{x=1}^{a+b-1} P_{a,b,c}(x,y) \left(x-\frac{a+b}{2}\right)^2
$$
and the vertical moment of inertia with respect to $S$
$$
\sum_{y=0}^{a+c-1} \sum_{x=1}^{a+b-1} P_{a,b,c}(x,y) 
\left(2\left(y-\frac{a+c-1}{2}\right)-\left(x-\frac{a+b}{2}\right)\right)^2
$$
of this distribution. (Note that 
$2(y-\frac{a+c-1}{2})-(x-\frac{a+b}{2})=0$ is the horizontal line 
in the coordinates of the oblique angled system which contains the lowest vertex of the horizontal 
lozenge in the centre of the hexagon.)
Our theorem is the following. 

\begin{theo}
\label{th:1}
Let $a$, $b$, $c$ be positive integers and let $P_{a,b,c}(x,y)$ denote the probability that an arbitrary chosen 
lozenge tiling of a hexagon with side lengths $a$, $b$, $c$, $a$, $b$, $c$ contains the 
horizontal lozenge with lowest vertex $(x,y)$. Then the horizontal moment of inertia 
with respect to $S$ is equal to
\begin{equation}
\label{horz}
\sum_{x=1}^{a+b-1} \sum_{y=0}^{a+c-1} P_{a,b,c}(x,y) \left(x-\frac{a+b}{2}\right)^2 = 
\frac{1}{12} a b (a^2 + b^2 -2) 
\end{equation}
and the vertical moment of inertia with respect to $S$ is equal to 
\begin{multline}
\label{vert}
\sum_{y=0}^{a+c-1} \sum_{x=1}^{a+b-1} P_{a,b,c}(x,y) 
\left(2\left(y-\frac{a+c-1}{2}\right)-\left(x-\frac{a+b}{2}\right)\right)^2 =  \\
= \frac{1}{12}  a b (a^2 + b^2 -2 + 4 c^2 + 4 a c + 4 b c).
\end{multline}
\end{theo}

In fact Propp \cite{propp} has already noticed that the computation of 
the horizontal moment of inertia is easy and states the formula for 
the case $a=b=c=n$. Furthermore, he concludes that the vertical moment 
of inertia is more difficult to compute for the first few values 
for the case $a=b=c=n$ do not seem to predict the formula to be 
a polynomial of degree $4$ as it is the case for the horizontal moment of 
inertia. However, this conclusion is based on a miscalculation. The vertical moment of inertia for the 
case $a=b=c=2$ is equal to $18$ (not to $20$) and thus we obtain the polynomial $7 n^4 / 6 - n^2 /6$ for the vertical moment of inertia if $a=b=c=n$. (See \eqref{vert}.)  
Nevertheless the computation of the vertical moment of inertia seems to be more involved compared to the 
computation of the horizontal moment of inertia.

\medskip

The following section is devoted to the proof of Theorem~\ref{th:1}. Our combinatorial proof is based on the 
correspondance between lozenge tilings of a hexagon with side lengths $a$, $b$, $c$, $a$, $b$, $c$ and plane partitions in an $a \times b \times c$ box, i.e.
plane partitions of shape $b^a$ and with entries between $0$ and $c$. Thus we can avoid to use an explicit 
expression for the probability $P_{a,b,c}(x,y)$.
In Lemma~\ref{2} we observe that we are in fact able to  compute the inner 
sum of \eqref{horz}. In order to obtain \eqref{vert}, we split up the left hand side of \eqref{vert} into three double sums and again we 
compute their inner sums after possibly interchanging the summation order (Lemma~\ref{1} - \ref{3}). In Lemma~\ref{4} we demonstrate that 
Lemma~\ref{3} and thus the computation of the vertical moment of inertia is trivial if we assume
$a=b$.

\section{Five lemmas and the proof of the Theorem}

\begin{lem}
\label{1}
Let $a$, $b$, $c$ be positive integers and $0 \le y_0 \le a+c-1$. Furthermore let $P_{a,b,c}(x,y)$ denote the probability that an arbitrary chosen 
lozenge tiling of the hexagon with side lengths $a$, $b$, $c$, $a$, $b$, $c$ contains the horizontal lozenge with lowest vertex $(x,y)$ in the 
oblique angled coordinate system. Then
$$\sum_{x=1}^{a+b-1} P_{a,b,c}(x,y_0) = \frac{a b}{a+c}.$$
\end{lem}

\smallskip

 The main ingredient for the proof is the symmetry of the Schur function.

\smallskip

{\it Proof.}
For a fixed $0 \le y_0 \le a+c-1$ the sum in question is just the expected value for 
the number of horizontal lozenges with its lowest vertex on the line $y = y_0$ in the oblique angled 
coordinate system. First we show that the sum does not depend on  $y_0$.
For a fixed plane partition in an $a \times b \times c$ box,
let $N_i(y_0)$ denote the multiplicity 
of the entry $y_0$ in the $i$th row of the plane partition. Then the number of horizontal lozenges 
with its lowest vertex on the line $y=y_0$ in a given lozenge tiling is equal to
$$
\sum_{i=1}^a N_i(y_0-a+i)
$$
in the corresponding plane partition. If we add $a-i$ to every entry in the $i$th row ($1 \le i \le a$) of the plane partition
we obtain a plane partition with strictly decreasing columns.
Then the sum above is just the number of $y_0$'s in this 
plane partition. The content $(\mu_i)_{i \ge 0}$ of a plane partition is the sequence with 
$$ \mu_i = \mbox{number of  i's in the plane partition}.$$ 
Fix a sequence $(\mu_i)_{i \ge 0}$ and an integer $j \ge 0$. Let 
$\nu_i = \mu_i$ if $i \not= j, j+1$, $\nu_j = \mu_{j+1}$ and 
$\nu_{j+1}=\mu_j$. 
By the bijection  in  \cite[page 152]{sagan} the number of plane partitions 
with decreasing rows,  strictly decreasing columns and  content 
$(\mu_i)_{i \ge 0}$ is equal to the number of such plane partitions with 
content  $(\nu_i)_{i \ge 0}$. Thus the expected value in question
is independent of $y_0$.
Now the assertion follows since the total number of horizontal 
lozenges is equal to $a b$, i.e.
$$\sum_{y=0}^{a+c-1} \sum_{x=1}^{a+b-1} P_{a,b,c}(x,y) = a b.$$
\qed

\bigskip

In our next lemma we compute the inner sum of \eqref{horz}.

\begin{lem}
\label{2}
Let $a$, $b$, $c$ be positive integers and $1 \le x_0 \le a+b-1$. Furthermore let 
$P_{a,b,c}(x_0,y_0)$ denote the probability that an arbitrary chosen lozenge tiling
of the hexagon with side lengths $a$, $b$, $c$, $a$, $b$, $c$ contains the 
horizontal lozenge with lowest vertex $(x,y)$ in the oblique angled coordinate system. Then
$$ \sum_{y=0}^{a+c-1} P_{a,b,c}(x_0,y) = 
\begin{cases}
\begin{array}{ll}
x_0  &  1 \le x_0 \le \min(a,b) \\
\min(a,b) & \min(a,b) \le x_0 \le \max(a,b) \\
a+b-x_0 & \max(a,b) \le x_0 \le a+b-1
\end{array}
\end{cases}. $$
\end{lem}

{\it Proof of Lemma~\ref{2}:}
For a fixed $1 \le x_0 \le a+b-1$ the sum in question is just the expected 
value for the number of horizontal lozenges with its lowest vertex on the 
vertical line $x=x_0$. But this number does not depend on the lozenge tiling. \qed

\begin{lem} 
\label{3}
Let $a$, $b$, $c$ be positive integers and $1 \le x_0 \le a+b-1$. Furthermore let 
$P_{a,b,c}(x,y)$ denote the probability that an arbitrary chosen lozenge tiling
of the hexagon with side lengths $a$, $b$, $c$, $a$, $b$, $c$ contains the 
horizontal lozenge with lowest vertex $(x,y)$ in the oblique angled coordinate system.
Then  
$$ \sum_{y=0}^{a+c-1} P_{a,b,c}(x_0,y) \left( y - \frac{a+c-1}{2} \right) =
\begin{cases}
\begin{array}{ll}
\frac{(-a^2- a b -a c + b c + a x_0 + b x_0) x_0}{2 (a+b)} & 1 \le x_0 \le a \\
\frac{a c (a + b - 2 x_0)}{2 (a+b)} & a \le x_0 \le b \\
\frac{(-b^2 - a b - b c + a c + a x_0 + b x_0) (a + b -x_0)}{2 (a+b)} & b \le x_0 \le a+b-1 
\end{array}
\end{cases}
$$
if $a \le b$ and 
$$ \sum_{y=0}^{a+c-1} P_{a,b,c}(x_0,y) \left( y - \frac{a+c-1}{2} \right) =
\begin{cases}
\begin{array}{ll}
\frac{(-a^2- a b -a c + b c + a x_0 + b x_0) x_0}{2 (a+b)} & 1 \le x_0 \le b \\
\frac{-b ( a+ b + c ) (a + b - 2 x_0)}{2 (a+b)} & b \le x_0 \le a  \\
\frac{(-b^2 - a b - b c + a c + a x_0 + b x_0) (a + b -x_0)}{2 (a+b)} & a \le x_0 \le a+b-1 
\end{array}
\end{cases}
$$
if $a \ge b$.
\end{lem}

\bigskip

For the proof of Lemma~\ref{3} we need the following definition of an $(n;k_1,k_2, \dots, k_a)$-array and the following lemma.
Roughly speaking an  $(n;k_1,k_2, \dots, k_a)$-array is the bottom-left part of a plane partition $T$ in an $a \times b \times c$-box 
with $T(i,i+n)=k_i$, which is dissected along the set of cells $\{(i,i+n),1 \le i \le a \}$. 

\begin{figure}
\setlength{\unitlength}{1cm}
\scalebox{0.35}{\includegraphics{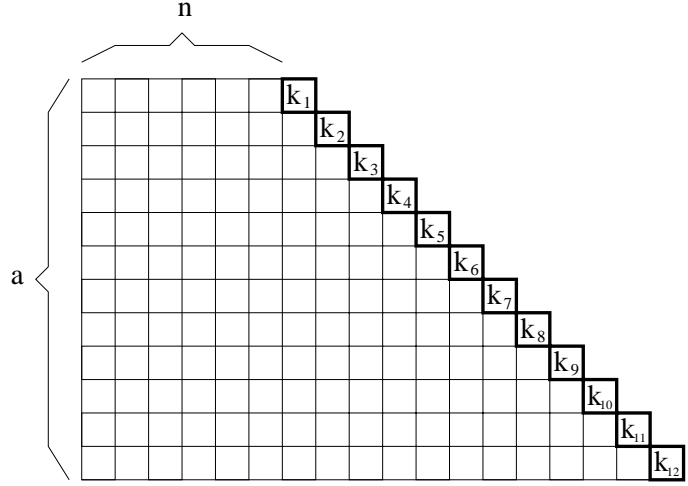}}
\caption{$F_{a,n}$ for $a=12$ and $n=6$.}
\label{fan}
\end{figure}

\begin{defi}
For every two positive integers $a$ and $n$ define the following set of cells:
$$
F_{a,n} = \{ (i,j) | 1 \le i \le a, 1 \le j \le i + n \}
$$
(See Figure~\ref{fan}.)
Let $(k_1,k_2,\dots, k_a)$ be a decreasing sequence of positive integers smaller or equal to $c$. Then an 
$(n;k_1,k_2,\dots,k_a)$-array is an assignment $T$ of integers to the cells in 
$F_{a,n}$ such that
\begin{enumerate}

\item[(i)]  $T(i,i+n)= k_i$ for $1\le i \le a$, 

\item[(ii)] $T(i,j) \le c$ for all $(i,j) \in F_{a,n}$ and

\item[(iii)] rows and columns are decreasing.
\end{enumerate}
The norm of an $(n;k_1,k_2,\dots,k_a)$-array $T$ is defined as 
$$
n(T) = \sum_{(i,j) \in F_{a,n}} T(i,j).
$$
Let $A_n(k_1,k_2,\dots,k_a)$ denote the number of $(n;k_1, k_2, \dots, k_a)$-arrays and 
$$
S_n(k_1,k_2,\dots,k_a) = \sum_{T}  \sum_{(i,j) \in F_{a,n}, \atop i +n \not= j} T(i,j),
$$
where the outer sum is taken over all $(n;k_1,k_2,\dots,k_a)$-arrays $T$.
\end{defi}

\begin{lem} 
\label{explizit}
Let $a$ and $n$ be positive integers and let $(k_1,k_2,\dots,k_a)$ be a decreasing sequence of positive integers smaller 
or equal to $c$. Then
$$
\frac{S_n(k_1,k_2,\dots,k_a)}{A_n(k_1,k_2,\dots,k_a)} = \frac{1}{2} \left( n \, a \, c + (a+n-1) \sum\limits_{i=1}^a k_i \right).
$$ 
\end{lem}

{\it Proof:} First we describe a bijection between $(n;k_1,k_2,\dots,k_a)$-assignments and 
semistandard tableaux of shape $(c-k_a,c-k_{a-1},\dots,c-k_1)$ with entries between $1$ and $a+n$.

Let $T$ be an arbitrary  $(n;k_1,k_2,\dots,k_a)$-assignment. In order to obtain the corresponding 
semistandard tableau replace every entry $T(i,j)$ of $T$ with  $c - T(i,j)$. We obtain an assignment with 
increasing rows and columns. Every row of this assignment is a partition, where 
the parts are written in reverse order. We conjugate these partitions and write them right justified and in increasing  order
among each other. Clearly the rows are weakly increasing but the columns are stricly increasing. If we rotate this assignment by $180^\circ$ we obtain a reverse semistandard 
tableau of shape $(c-k_{a},c-k_{a-1},\dots,c-k_1)$ with 
entries between $1$ and $a+n$.  Finally we replace every entry $e$ of this reverse semistandard tableau by $a+n+1-e$ and obtain the desired semistandard tableau $T'$. 
This procedure is obviously reversible.  

The norms of the two corresponding objects are related as follows
$$
n(T) = c \left( a \, n + \frac{ (a+1) \, a}{2} \right) - (a + n +1) \left( a \, c - \sum_{i=1}^a k_i \right) + n(T'),
$$
where $n(T')$ is the sum of the parts of $T'$.

Note that the statement in the lemma is equivalent to
$$
\frac{\sum\limits_{T} n(T)}{\sum\limits_{T} 1} = \frac{1}{2}  \left( n \, a \, c + (a+n+1) \sum\limits_{i=1}^a k_i \right),
$$
where the sums are taken over all $(n;k_1,k_2,\dots,k_a)$-arrays $T$. Thus we have to show that 
$$
\frac{\sum\limits_{T'} n(T')}{\sum\limits_{T'} 1} = \frac{1}{2} (a \, c - \sum_{i=1}^a k_i) (a+n+1),
$$
where the sums are taken over all semistandard tableaux $T'$ of shape $(c-k_a,c-k_{a-1},\dots,c-k_1)$ with entries between $1$ and $a+n$ or, 
equivalently, that 
$$
\frac{\sum\limits_{T'} n(T')}{\sum\limits_{T'} 1} = \frac{1}{2} \left( \sum_{i=1}^{r} \lambda_i \right) (a+1),
$$
where the sums are taken over all semistandard tableaux $T'$ of shape $(\lambda_1, \lambda_2, \dots, \lambda_r)$ with entries between $1$ and 
$a$.

By Stanley's hook-content formula \cite[Theorem~15.3]{stan} we 
have
\begin{equation}
\label{stan}
\sum_{T'} q^{n(T')} = q^{\sum\limits_{i=1}^r i \lambda_i} \prod_{\rho \in \lambda} \frac{1 - q^{a + c_{\rho}}}{1-q^{h_{\rho}}},
\end{equation}
where the sum is taken over all semistandard tableaux $T'$ of shape $\lambda:=(\lambda_1, \lambda_2, \dots, \lambda_r)$ with entries between $1$ and
$a$ and the product is taken over all cells $\rho$ in the Ferrer diagramm of shape $\lambda$. Furthermore $h_{\rho}$ denotes the 
hooklength and $c_{\rho}$ the content of cell $\rho$. We observe that 
$$
\frac{\sum\limits_{T'} n(T')}{\sum\limits_{T'} 1} = \lim_{q \to 1}  \frac{ \left( \sum\limits_{T'} q^{n(T')} \right)' }{\sum\limits_{T'} q^{n(T')}}.
$$
With the help of \eqref{stan} we compute the derivative of $\sum\limits_{T'} q^{n(T')}$.
\begin{eqnarray*}
\left (\sum_{T'} q^{n(T')} \right)' 
&=& \left( \sum_{i=1}^r i \lambda_i \right) q^{(\sum\limits_{i=1}^r i \lambda_i) - 1} \prod_{\rho \in \lambda} \frac{1 - q^{a + c_{\rho}}}{1-q^{h_{\rho}}} 
+ q^{\sum\limits_{i=1}^r i \lambda_i } \\
&& \times \sum_{\rho \in \lambda} \frac{-(a+c_{\rho}) (q^{a+c_{\rho} - 1}) ( 1 - q^{h_{\rho}}) + ( 1 - q^{a+c_{\rho}}) ( h_{\rho} q^{h_\rho - 1} )}
                              { (1 - q^{h_{\rho}}) ( 1 - q^{a + c_{\rho}}) } \\
&& \times \prod_{\rho' \in \lambda} \frac{1 - q^{a + c_{\rho'}}}{1-q^{h_{\rho'}}}
\end{eqnarray*}
Thus 
\begin{eqnarray*}
\lim_{q \to 1}  \frac{ (\sum\limits_{T'} q^{n(T')})' }{\sum\limits_{T'} q^{n(T')}} &=&
\sum_{i=1}^r i \lambda_i + \frac{1}{2} \sum_{\rho \in \lambda} \left( a + c_{\rho} - h_{\rho} \right) \\
&=& \sum_{i=1}^r i \lambda_i + \frac{a}{2} \sum_{i=1}^r \lambda_i + \frac{1}{2} \sum_{j=1}^{\lambda_1} j \lambda'_j - 
                                                               \frac{1}{2} \sum_{i=1}^r i \lambda_i    \\
&& -
                                                               \frac{1}{2} \left( \sum_{i=1}^r \frac{( \lambda_i + 1 ) \lambda_i}{2} + 
                                                                                  \sum_{j=1}^{\lambda_1}  \frac{ \lambda'_j  (\lambda'_j - 1)}{2}  \right) \\
&=& \frac{(a+1)}{2} \sum_{i=1}^r \lambda_i,
\end{eqnarray*}
for 
$$
\sum_{i=1}^r i \lambda_i = \sum_{j=1}^{\lambda_1} \frac{(\lambda'_j + 1) \lambda'_j}{2}  
$$
are two ways to express the norm of the tableau of shape $\lambda$ and with constant entry $i$ in every cell in the $i$th row and 
$$
\sum_{j=1}^{\lambda_1} j \lambda'_j = \sum_{i=1}^r \frac{(\lambda_i + 1) \lambda_i}{2}
$$
are two way to express the norm of the tableau of shape $\lambda$ and with constant entry $j$ in every cell in the $j$th column.
This concludes the proof of the lemma. \qed

\medskip

In the following $E(i,j)$ denotes the expected values of the entry in cell $(i,j)$ of a plane partition in an $a \times b \times c$ box.

\smallskip

{\it Proof of Lemma~\ref{3}.} If $a \le b$ we have
$$
\sum_{y=0}^{a+c-1} P_{a,b,c} (x,y) \,\, y =
\begin{cases}
\begin{array} {ll}
\sum\limits_{i=1+a-x}^a \left( E(i,x-a+i) + a - i \right) & 1 \le x \le a \\
\sum\limits_{i=1}^a \left( E(i,x-a+i) + a -i \right) & a \le x \le b \\
\sum\limits_{i=1}^{a+b-x} \left( E(i,x-a+i) + a - i \right) & b \le x \le a+b-1 
\end{array}
\end{cases}
$$
and if $b \le a$ we have
$$
\sum_{y=0}^{a+c-1} P_{a,b,c} (x,y) \,\, y =
\begin{cases}
\begin{array} {ll}
\sum\limits_{i=1+a-x}^{a} \left( E(i,x-a+i) + a - i \right) & 1 \le x \le b \\
\sum\limits_{i=1+a-x}^{a+b-x} \left( E(i,x-a+i) + a -i \right) & b \le x \le a \\
\sum\limits_{i=1}^{a+b-x} \left( E(i,x-a+i) + a - i \right) & a \le x \le a+b-1 
\end{array}
\end{cases}.
$$
Thus, by Lemma~\ref{2}, we have to show that 
$$
\sum_{i=1+a-x}^a E(i,x-a+i) = \frac{b \, c \, x}{a + b} 
$$
for $1 \le x \le \min(a,b)$
and 
$$
\sum_{i=1}^{a+b-x} E(i,x-a+i) =  \frac{a \, c \, (a + b - x)}{a+ b}
$$
for $\max(a,b) \le x \le a + b -1$.
Furthermore, if $a \le b$ we have to show that
$$
\sum_{i=1}^a E(i,x-a+i) = \frac{a \, c \, (a + b - x)}{a+ b} 
$$
for $a \le x \le b$ and if $b \le a$ we have to show that 
$$
\sum_{i=1+a-x}^{a+b-x}  E(i,x-a+i) = \frac{b \, c \, x}{a + b} 
$$
for $b \le x \le a$.

In the following we will only consider the case $a \le b$ for the other case is similar. Let
$$
D(x) =
\begin{cases}
\begin{array} {ll}
\sum\limits_{i=1+a-x}^a E(i,x-a+i) & 1 \le x \le a  \\
\sum\limits_{i=1}^a E(i,x-a+i) & a \le x \le b \\
\sum\limits_{i=1}^{a+b-x} E(i,x-a+i) & b \le x \le a+b-1
\end{array}
\end{cases}.
$$

\begin{figure}
\setlength{\unitlength}{1cm}
\scalebox{0.35}{\includegraphics{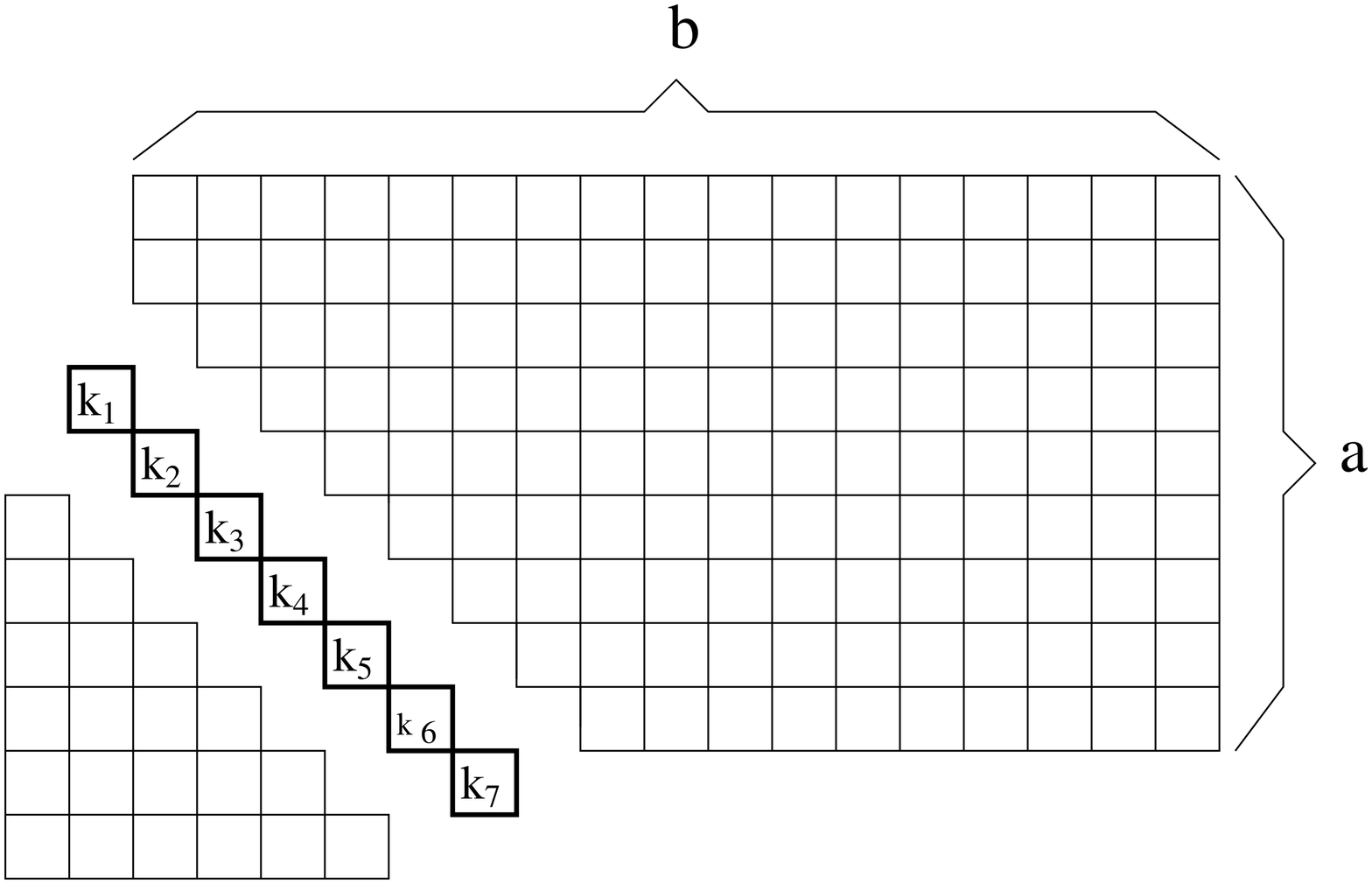}}
\caption{}
\label{diag1}
\end{figure}

Let $1 \le x \le a$  and $G_x$ be the following set of cells
$$
G_x = \{(i,j) | 1 \le i \le a, 1 \le j \le b\} \setminus \{ (i, y-a + i) | 1  \le y <  x , 1+ a -y \le i \le a \}.
$$
(See Figure~\ref{diag1}.) An $x$-partial plane partition in an $a \times b \times c$ box is an assignment of the cells in $G_x$ with integers 
between  $0$ and $c$ such that rows and columns are decreasing. Let $(k_1,k_2,\dots,k_x)$ be a decreasing sequence of  
integers in $\{0,1,\dots,c\}$ and  let $U(k_1,k_2,\dots,k_x)$ be the number 
of $x$-partial plane partitions $T$ with $T(i,x-a+i) = k_{x-a+i}$ for all $i$. Furthermore let $N(a,b,c)$ denote the number of 
plane partitions in an $a \times b \times c$ box.
Then
$$
D(x) = \frac{1}{N(a,b,c)} \sum_{(k_1,k_2,\dots,k_x)} A_0(k_1,k_2,\dots,k_x) \,\, U(k_1,k_2,\dots,k_x)\sum_{i=1}^x k_i 
$$
for $1 \le x \le a$.
(See Figure~\ref{diag1}.)

Thus 
\begin{multline*}
\frac{1}{x} D(x) =  \frac{1}{N(a,b,c)} \sum_{(k_1,k_2,\dots,k_x)} \frac{2}{x (x-1)} S_0(k_1,k_2,\dots,k_x) \,\, U(k_1,k_2,\dots,k_x) \\
= \frac{1}{\binom{x}{2}} \sum_{y=1}^{x-1} D(y)  
\end{multline*}
for $1 \le x \le a$ by Lemma~\ref{explizit}.  By induction we obtain that the fraction
$\frac{D(x)}{x} $
is independent of $x$ for $1 \le x \le a$. By symmetry 
$\frac{D(x)}{a+b-x}$
is independent of $x$ for $b \le x \le a+b-1$.
Thus 
\begin{equation}
\label{d1}
D(x) = x \, D(1)
\end{equation}
for $1 \le x \le a$ and 
\begin{equation}
\label{d2}
D(x) = (a+b-x) \, D(a+b-1)
\end{equation}
for $b \le x \le a+b-1$.

\begin{figure}
\setlength{\unitlength}{1cm}
\scalebox{0.35}{\includegraphics{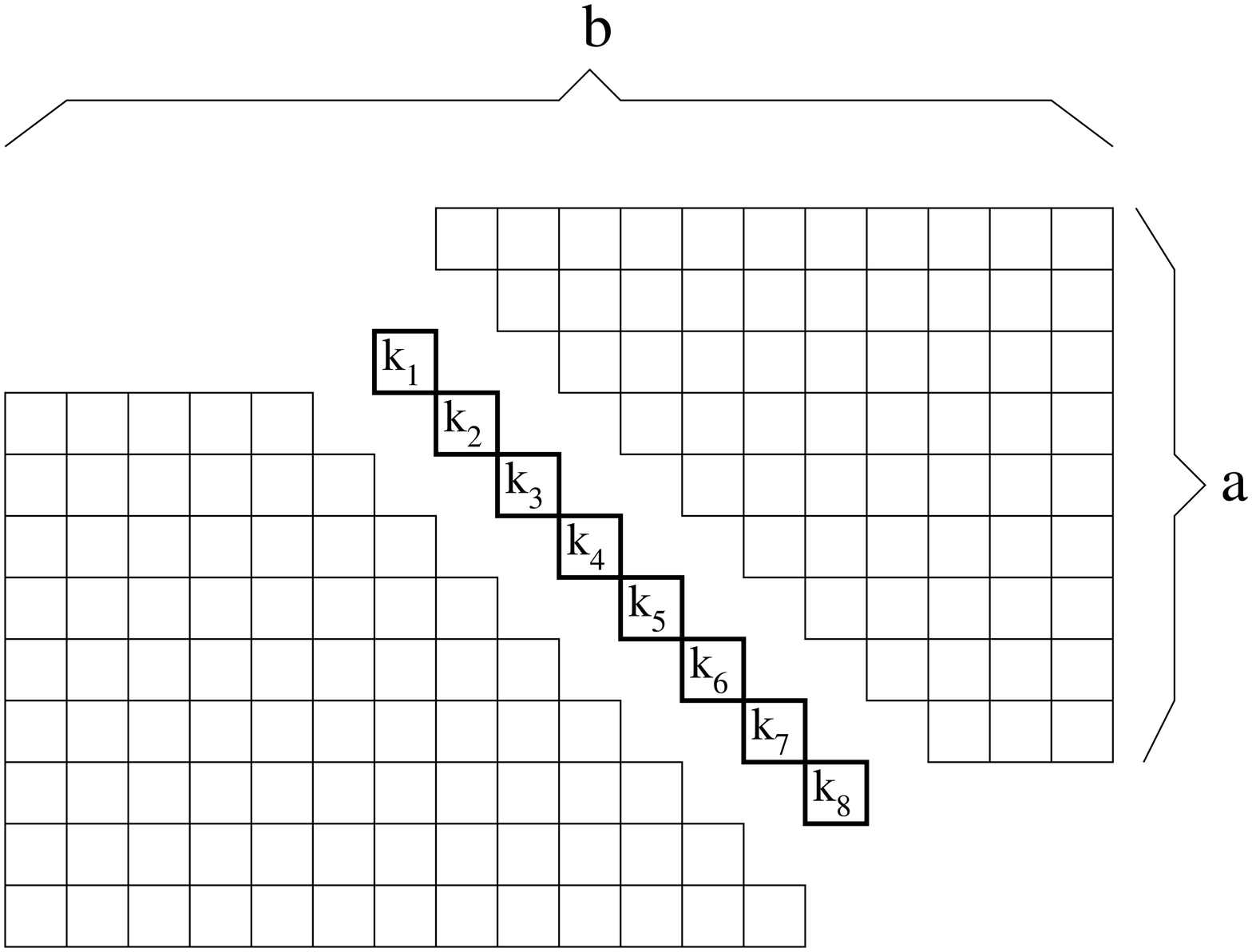}}
\caption{}
\label{diag2}
\end{figure}

Let $a \le x \le b$ and $G_x$ be the following set of cells.
\begin{multline*}
G_x = \{(i,j) | 1  \le i  \le a, 1 \le j \le b\} \setminus \\ \{(i, y-a+i) | (1 \le y < a, 1+a-y \le i \le a) \vee (a \le y <  x, 1 \le i \le a)\}
\end{multline*}
(See Figure~\ref{diag2}.)
An $x$-partial plane partition is defined in the same way as above. Let $(k_1,k_2,\dots,k_a)$ be a decreasing sequence of integers
in $\{0,1,\dots,c\}$ and $U(k_1,k_2,\dots,k_a)$ the number of $x$-partial plane partitions $T$ with $T(i,x-a+i)=k_i$ for $1 \le i \le a$.

Then, by Lemma~\ref{explizit},  we have  (see Figure~\ref{diag2})
\begin{eqnarray*}
D(x)
&=& \frac{1}{N(a,b,c)} \sum_{(k_1,k_2,\dots,k_x)} A_{x-a}(k_1,k_2,\dots,k_a) \, U(k_1,k_2,\dots,k_a) \sum_{i=1}^a k_i \\
&=& \frac{1}{N(a,b,c)} 
 \sum_{(k_1,k_2,\dots,k_x)} \frac{2 S_{x-a}(k_1,k_2,\dots,k_a) - A_{x-a}(k_1,k_2,\dots,k_a)\, (x-a) \, a \, c}{(x-1)} \\
&& \qquad \times  U(k_1,k_2,\dots,k_a) \\
&=& \frac{2}{x-1} \sum_{y=1}^{x-1} D(y)   - \frac{(x-a) \, a \, c}{(x-1)} \\
&=& \frac{2}{x-1} \sum_{y=a+1}^{x-1} D(y) + \frac{D(1) \, a \, (a+1)}{x-1} - \frac{(x-a) \, a \,c}{x-1}
\end{eqnarray*} 
and therefore, by induction,
\begin{equation}
\label{d3}
D(x) = x \, D(1) - (x-a) \, c
\end{equation}
for $a \le x \le b$.
On one hand we have 
$$
D(b)=b \, D(1) - (b-a) \, c
$$
by \eqref{d3}.
But on the other we have
$$
D(b) = a \, D(a+b-1)
$$
by \eqref{d2}.
Furthermore we have
$$
D(1) + D(a+b-1) = E(a,1) + E(1,b) =c 
$$
by the involution on the set of plane partitions in an  $a \times b \times c$ box with 
$T(i,j) = c - T(a+1-i,a+1-j)$, where $T(i,j)$ denotes the entry of the cell $(i,j)$.
Thus 
$$
D(1)=\frac{b c}{a+b}
$$
and 
$$
D(a+b-1) = \frac{a c}{a+b}
$$
and the assertion follows by \eqref{d1}, \eqref{d2} and \eqref{d3}.

\bigskip

The proof of the following lemma shows that Lemma~\ref{3} is nearly obvious if we assume $a=b$.

\begin{lem} 
\label{4}
Let $m$, $n$ be positive integers and $1 \le x_0 \le 2n-1$. Then
$$ \sum_{y=0}^{n+m-1} P_{n,n,m}(x_0,y) \left( y - \frac{n+m-1}{2} \right) =
\begin{cases}
\begin{array}{ll}
\frac{1}{2} (x_0-n) x_0 & 1 \le x_0 \le n \\
\frac{1}{2} (x_0-n) (2n -x_0) & n \le x_0 \le 2n 
\end{array}
\end{cases}.
$$
\end{lem}

\bigskip

\bigskip

{\it Proof of Lemma~\ref{4}:}
By Lemma~\ref{2} we have
$$
\sum_{y=0}^{n+m-1} P_{n,n,m}(x_0,y)  \left( - \frac{n+m-1}{2} \right) =
\begin{cases}
\begin{array} {ll}
\frac{1-n-m}{2} \,\, x_0 & 1 \le x_0 \le n \\
\frac{1-n-m}{2} \,\, (2n - x_0) &  n \le x_0 \le 2n - 1 
\end{array}
\end{cases}.
$$
Thus we have to show that  
$$
\sum_{y=0}^{n+m-1} P_{n,n,m}(x_0,y) \,\, y =
\begin{cases}
\begin{array} {ll}
\frac{x_0+m-1}{2} \,\, x_0 & 1 \le x_0 \le n \\
\frac{x_0+m-1}{2}\,\, (2n - x_0) &  n \le x_0 \le 2n - 1 
\end{array}
\end{cases}.
$$
This sum is equal to
$$
\sum_{y=0}^{n+m-1} P_{n,n,m} (x_0,y) \, \, y  =
\begin{cases}
\begin{array} {ll}
\sum\limits_{i=n+1-x_0}^n \left( E(i,x_0-n+i) + n-i \right) & 1 \le x_0 \le n \\
\sum\limits_{i=1}^{2n - x_0} \left( E(i,x_0-n+i) + n-i  \right) & n \le x_0 \le 2n - 1 
\end{array}
\end{cases},
$$
where $E(i,j)$ is the expected value of the entry in cell $(i,j)$
in a plane partition in an $n \times n \times m$ box.
Therefore we have show that
$$
\sum_{i=n+1-x_0}^n E(i,x_0-n+i)  = \frac{1}{2} \, x_0 \, m  \quad \mbox{if} \quad  1 \le x_0 \le n 
$$
and
$$
\sum_{i=1}^{2n - x_0} n-i = \frac{1}{2} \, (2n - x_0) \, m  \quad \mbox{if} \quad   n \le x_0 \le 2n - 1. 
$$
But this is obvious since $E(i,j) + E(n+1-j,n+1-i) = m$ by the involution 
on the set of plane partitions in an $n \times n \times m $ box 
with $T(i,j) = m - T(n+1-j,n+1-i)$, 
where $T(i,j)$ denotes the entry of cell $(i,j)$. \qed

\bigskip

Finally we combine Lemma~\ref{1} -- \ref{3} in order to prove Theorem~\ref{th:1}. 

\medskip

{\it Proof of Theorem~\ref{th:1}:}
In order to calculate the horizontal moment of inertia we use Lemma~\ref{2}.
\begin{multline*}
\sum_{x=1}^{a+b-1} \left( \sum_{y=0}^{a+c-1} P_{a,b,c} (x,y) \right) \left( x - \frac{a+b}{2} \right)^2  \\
= \sum_{x=1}^{\min(a,b)} x \left( x - \frac{a+b}{2} \right)^2 + \sum_{x=\min(a,b)+1}^{\max(a,b)-1} \min(a,b) 
\left( x - \frac{a+b}{2} \right)^2 \\  + \sum_{x=\max(a,b)}^{a+b-1} (a+b-x) \left( x  - \frac{a+b}{2} \right)^2 \\
= \frac{1}{12} a b  (a^2 + b^2 -2) 
\end{multline*}

The vertical moment of inertia splits into 3 sums.
\begin{multline}
\label{4sum}
\sum_{y=0}^{a+c-1} \sum_{x=1}^{a+b-1} P_{a,b,c}(x,y) 
\left(2\left(y-\frac{a+c-1}{2}\right)-\left(x-\frac{a+b}{2}\right)\right)^2  \\
= 4 \sum_{y=0}^{a+c-1} \left( \sum_{x=1}^{a+b-1} P_{a,b,c} (x,y) \right) \left( y - \frac{a+c-1}{2} \right) ^2  \\
 - 4 \sum_{x=1}^{a+b-1} \left( \sum_{y=0}^{a+c-1} P_{a,b,c}(x,y) \left( y - \frac{a+c-1}{2} \right) \right) \left( x - \frac{a+b}{2} \right) \\
 + \sum_{x=1}^{a+b-1} \left( \sum_{y=0}^{a+c-1} P_{a,b,c}(x,y) \right) \left( x - \frac{a+b}{2} \right)^2 
\end{multline}
The last sum is equal to the horizontal moment of inertia.

In order to compute the first sum in \eqref{4sum} we use Lemma~\ref{1}.
\begin{multline*}
4 \sum_{y=0}^{a+c-1} \left( \sum_{x=1}^{a+b-1} P_{a,b,c} (x,y) \right) \left( y - \frac{a+c-1}{2} \right) ^2 \\
= 4 \sum_{y=0}^{a+c-1} \frac{a b}{a+c} \left( y - \frac{a+c-1}{2} \right)^2  \\
= \frac{1}{3}  a b ( a + c - 1) (a + c + 1)
\end{multline*}

For the second sum in \eqref{4sum} we use Lemma~\ref{3}.
\begin{multline*}
4 \sum_{x=1}^{a+b-1} \left( \sum_{y=0}^{a+c-1} P_{a,b,c}(x,y) \left( y - \frac{a+c-1}{2} \right) \right) \left( x - \frac{a+b}{2} \right) \\
= \frac{1}{3} a b (-1 + a^2 + a c - b c )
\end{multline*}

Thus 
\begin{multline*}
\sum_{y=0}^{a+c-1} \sum_{x=1}^{a+b-1} P_{a,b,c}(x,y) 
\left(2\left(y-\frac{a+c-1}{2}\right)-\left(x-\frac{a+b}{2}\right)\right)^2  \\
= \frac{1}{3} a b c (a + b + c )  + 
\sum_{x=1}^{a+b-1} \left( \sum_{y=0}^{a+c-1} P_{a,b,c} (x,y) \right) \left( x - \frac{a+b}{2} \right)^2
\end{multline*}
and the assertion follows.
\qed


\begin{thebibliography} {2}

\bibitem{propp}
J. Propp, Enumeration of Matchings: Problems and Progress, {\it New Perspectives in 
Geometric Combinatorics}, Edited by L. Billera, A. Bj\"oner, C. Greene, R. Simeon, 
R. Stanley, Mathematical Sciences Research Institute Publications {\bf 38}, Cambridge 
University Press, 1999.

\bibitem{sagan}
B. Sagan, The Symmetric Group: Representations, Combinatorial Algorithms and Symmetric Functions,
Wadsworth and  Brooks/Cole, Mathematics Series, 1991. 

\bibitem{stan}
R. P. Stanley, {\it Theory and applications of plane partitions: Part 2}, Stud. Appl. Math. {\bf 50} (1971), 259 - 279.

\end{thebibliography}
\end{document}